%Feb 28, 2002
\input amstex
\documentstyle{amsppt}
\magnification=1200
\hsize=13.8cm
\vsize=19 cm
\catcode`\@=11
\def\NoLogo{\let\logo@\empty}
\catcode`\@=\active
\NoLogo

\def \b {\beta}
\def\i{\sqrt{-1}}
\def\Ric{\text{Ric}}
\def\lf{\left}
\def\ri{\right}
\def\bbar{\bar \beta}
\def\a{\alpha}
\def\g{\gamma}
\def\p{\partial}
\def\delbar{\bar\delta}
\def\ddbar{\partial\bar\partial}

\def\C{\Bbb C}
\def\R{\Bbb R}

\def\vp{\varphi}
\def\tD{\tilde \Delta}

\def\bb{{\bar\beta}}
\def\abb{{\alpha\bar\beta}}

\def\i{\sqrt {-1}}
\def\tD{\widetilde \Delta}
\def\tn{\widetilde \nabla}
\def \D {\Delta}

\def\aint{\frac{\ \ }{\ \ }{\hskip -0.4cm}\int}
\documentstyle{amsppt}
\magnification=1200
\hsize=13.8cm

\leftheadtext{Lei Ni  and Luen-fai Tam}
\rightheadtext{K\"ahler-Ricci flow, Poincar\'e-Lelong equation}
\topmatter

\title{K\"ahler-Ricci flow and the Poincar\'e-Lelong equation}\endtitle

\author{Lei Ni\footnotemark and  Luen-Fai Tam\footnotemark}\endauthor
\footnotetext"$^{1}$"{Research partially supported by NSF grant
DMS 0196405 and DMS-0203023, USA.}
\footnotetext"$^2$"{Research partially supported by   Earmarked Grant of
Hong Kong \#CUHK4217/99P.}

\address
Department of Mathematics, Stanford University, Stanford, CA 94305
\endaddress
\email{
lni\@math.stanford.edu}
\endemail

\address
Department of Mathematics, The Chinese University of Hong Kong,
Shatin, Hong Kong, China
\endaddress
\email{lftam\@math.cuhk.edu.hk}
\endemail

\affil
{
Stanford University\\
The Chinese University of Hong Kong
}
\endaffil

\date  September, 2001; revised in Feburary, 2002
\enddate

\endtopmatter

\document

\subheading{\S0 Introduction}

In \cite{M-S-Y}, Mok-Siu-Yau studied complete K\"ahler manifolds with
nonnegative holomorphic bisectional curvature by solving the
Poincar\'e-Lelong equation
$$
\i\ddbar u=\Ric\tag0.1
$$
where $\Ric$ is the Ricci form of the manifold. In \cite{M-S-Y},
 the authors solved (0.1) under the assumptions
that the manifold is of maximal volume growth and  the scalar curvature decays
quadratically. On the other hand, in a series of papers of W.-X. Shi \cite{Sh2-4},
K\"ahler-Ricci flow
$$
\frac{\p}{\p t}g_\abb=-R_\abb\tag0.2
$$
has been studied extensively and important applications were given.
In  [N1] and \cite{N-S-T}, the Poincar\'e-Lelong equation has been solved under
more general conditions than in \cite{M-S-Y}. The conditions in [N-S-T]
are more in line
with the conditions in \cite{Sh2-4}. Since a solution of (0.1) is a potential for
 the Ricci tensor, it is interesting to see if one can apply (0.1) to study
solutions of (0.2).

In this work, on the one hand we shall study the K\"ahler-Ricci flows by
using solutions of the Poincar\'e-Lelong equation. On the other hand, we will
also  refine some of the results in \cite{Sh3, C-Z, C-T-Z} and
give new applications.
The hinge between the equations (0.1) and (0.2)
is that by solving (0.1) one can then
construct a function $u(x,t)$ which satisfies the time-dependent heat equation
$(\frac{\p}{\p t}-\Delta )u(x,t)=0$
and the time-dependent Poincar\'e-Lelong equation $\i\ddbar u=\Ric_{g(t)}$
simultaneously. It then
can simplify the study of (0.2) quite a bit.
It also suggests some of the refined
estimates in the second part of this paper. We should point out here that the
simplification in this paper
 is that $|\nabla u|^2$ helps to obtain a sharp uniform curvature
estimates (Cf. Theorem 1.3),
which holds as an equality for the K\"ahler-Ricci soliton.
It is different from the compact case  as in [Co1], where
one  restricts the deformation of the metric within a fixed cohomology class
and can then appeal to Yau's
solution to the Monge-Amper\'e equation by reducing (0.2) to a single equation.

Let $(M^m,g_\abb(x))$ be a complete noncompact K\"ahler manifold with bounded and
nonnegative holomorphic bisectional curvature. Let $\Cal R_0$ be the scalar
curvature of $M$. In \cite{Sh3}, it was proved that (0.2) has long time solution
with initial metric $g_\abb(x)$ satisfying  the assumption that
$$
k(x,r)\le C(1+r)^{-\theta}\tag0.3
$$
for some constants $C$ and $\theta>0$ for all $x$ and $r$. Here $k(x,r)$ denotes
the average of $\Cal R_0$ on $B(x,r)$, the geodesic ball of radius $r$ with center
at $x$. The idea of the proof of the long time existence in \cite{Sh3} is to use the parabolic version of the
third derivative estimate for the Monge-Amper\'e equation together with a careful
estimate  of the volume element. The computation is rather tedious. In this work,
we will use the solution to
 (0.1) constructed in [N-S-T] (more precisely the uniform curvature
estimate (1.24) in Theorme 1.3) to  give an alternate (and much
 simpler, we believe) proof for the long time existence under the assumption that
$$
\int_0^\infty k(x,r)dr\le C\tag0.4
$$
for some $C$ independent of $x$. Our proof  uses a maximum principle
which is a generalization of that in \cite{K-L}, and an idea similar
to those in [Cw]. Our assumption here is
different from but somewhat stronger than Shi's (0.3). However it
has covered the
interesting cases in [Sh2-3], namely the cases $k(x,r)\le C(1+r)^{-1-\delta}$,
on which interesting geometric results could be obtained.
On the other hand
we also  can prove a long time existence result under a
more flexible condition.
Namely, we show that there exists long time solution to (0.2) if
$$
k(x,r)\le \epsilon(r)\tag0.5
$$
for {\it all} $x$ (with some fixed function  $\epsilon(r)$) with
$\epsilon(r)\to 0$ as $r\to \infty$.
Recently  in \cite{C-T-Z}, it is proved that if the complex dimension of $M$
is $m=2$ and $M$ has maximal volume growth,  then (0.2) has long time solution if
 (0.5) holds for {\it some} $x$ and for some function $\epsilon(r)$ which
 tends to zero as $r\to\infty$. The proof there is an indirect blow-up
argument. It also used
some special features in dimension 2, such as the Guass-Bonnet formula for the four
 dimensional Riemannian manifolds.
In order to prove the long time existence under the assumption (0.5),
we need a more precise estimate for the volume element
$F(x,t)=\log\lf[\det(g_\abb(x,t))/\det(g_\abb(x,0))\ri]$, where $g_\abb(x,t)$
 is the solution of (0.2). In fact, we prove the following results, see Theorem 2.1
 and Corollary 2.1:
\proclaim{Theorem} Suppose (0.2) has a solution on $M\times[0,T)$. Then we have
the following:

(a) There exists a constant $C>0$ depending
only on $m$ such that for $0<t<T$,
$$
-F(x_0,t)\ge C\int_0^{\sqrt t}sk(x_0,s)ds.
$$

(b) If in addition, $k(x,r)\le k(r)$ for some function $k(r)$ for all $x$, then
$$
-\frak m(t)\le C'\int_0^{R}sk(s)ds
$$
where $R^2=at(1-\frak m(t))$, $C$ and $a$ are constants depending only on $m$.
Here $\frak m(t)=\inf_{x\in M}F(x,t)$.
\endproclaim
From the two-sidedness of the above estimates on $F(x,t)$ one can see that
they are almost optimal.  By comparing with the previous
estimates obtained in [N-S-T] and [N2] for the Poisson equation and the
linear heat equation,
the refined estimates here are sharp in certain cases and fit into
the theory for the linear equation.
The above mentioned estimates will be proved by using,
the by-now standard estimates on the
heat kernels of Li-Yau in \cite{L-Y}. There is no need to construct special
exhaustion functions as in \cite{Sh2-3, C-Z, C-T-Z}.
As a consequence, a little  more general  gap theorem, than those
in \cite{C-Z},  is obtained, see Corollary 2.3. In particular, we show that
{\it any bounded  solution to the Poisson equation
$\Delta u={\Cal R}_0(x)$ is a constant, provided $M$ has
bounded nonnegative
bisectional curvature.}  In other words, if $M$ is nonflat,
$\Delta u={\Cal R}_0(x)$ has no bounded solution.
This  answers a question asked by R. Hamilton. Namely, solving
Poisson for ${\Cal R}_0(x)$ is  different from arbitrary $f(x)$ since one
can easily construct bounded solution to $\Delta u=f(x)$ for nonzero compact
supportted $f(x)$.
This is also related
 to
the gradient estimates of Chow in [Cw]. In [Y], it was proved that, on a
complete Riemannian manifold with nonnegative Ricci curvature, any
negative (positive) harmonic function is a constant. We prove that a similar
result holds for $\Delta u ={\Cal R}_0(x)$. Namely,
{\it $\Delta u ={\Cal R}_0(x)$ has no nonconstant negative
solution,  provided $M$ has
bounded nonnegative
bisectional curvature and (0.2) has long time solution.}

When $(M,g_{\abb}(x,0))$ has the
maximum
volume growth, using the estimates mentioned above  the results
in [C-Z] on the Steinness and the topology of $M$ can be refined.
Namely we show that
if $(M, g_{\abb}(x,0))$ is of maximum volume growth and
$\int_0^rsk(x, s)\, ds \le \phi(r)$ with $\phi(r)\to 0$ as $r\to \infty$,
$M$ is Stein and diffeomorphic to $\R^{2m}$ for $m\ge 3$,
homeomorphic to $\R^4$
for $m=2$.

Another application of the estimates of $F$ and  (0.1) is that one can prove
the preservation of the decay rate  of $\Cal R_0$ in a certain sense.
For example, we will prove in Theorem 2.3 that if $\int_0^rsk(x,s)ds\le
C\log(1+r)$ (or $C(1+r)$), where $k(x,r)$ is the average of the
scalar curvature at $t=0$,
then we still have $\int_0^rsk_t(x,s)ds\le C'\log(1+r)$
($C'(1+r)$, respectively), where $k_t(x,r)$ is
the average of the scalar curvature at time $t$. Note that the constant $C'$
is independent of $t$. This might be useful in analyzing the
singularity models obtained by the blow-up procedure as in [H3].

From the methods of proof of the estimates of $F$, we can show that, under a
 rather weak decay condition on $\Cal R_0$, the volume growth is preserved in
 the sense that for any $t>0$,
$$
\lim_{r\to\infty}\frac{V_t(o,r)}{V_0(o,r)}=1
$$
where $V_t(o,r)$ is the volume of the geodesic ball with center at $o$ and
radius $r$ with respect to $g_\abb(x,t)$. This generalizes the results of
\cite{H3, Sh2, C-Z, C-T-Z}.

In \cite{Sh2}, under the assumption that $\theta=2$ in (0.3) and that $M$ has
positive holomorphic bisectional curvature, Shi proved that the rescaled
metric $\widehat g_\abb(x,t)=g_\abb(x,t)/g_{v\bar v}(x_0,t)$ subconverges to
a flat K\"ahler metric on $M$, where $x_0$ is a fixed point and $v$ is a
fixed nonzero (1,0) vector at $x_0$. If $M$ has maximal volume growth and if
the limit metric is complete, then one can conclude that $M$ is biholomorphic
to $\Bbb C^m$. It is pointed out in \cite{C-Z} that from \cite{Sh2} it is
 unclear why the property of completeness is true. In Proposition 3.1, we
will prove that if the scalar curvature $\Cal R_0$ has pointwise  quadratic
decay, then the largest eigenvalue of the limit metric with respect to the
initial metric grows at least like $r_0^a(x)$ for some $a>0$, where $r_0(x)$
 is the distance function to a fixed point with respect to the initial metric.
 This is a consequence of  the
result that   volume elements of the rescaled metrics
converge to the solution of
the Poincar\'e-Lelong equation constructed in [N-S-T], see Theorem 3.1.
We believe that this new piece of information will be helpful in studying  the
completeness of the limiting metric.

Here is how we organize this paper. In \S1, we will give an alternate proof of long time existence for (0.2). In \S2, we will give more refined estimates for $F(x,t)$ together with some applications. In \S3, we will study the asymptotic behavior of $F(x,t)$.

We shall use the differential inequalities for K\"ahler-Ricci flow
of Cao \cite{Co2-3} from time
to time, which is also called Harnack inequality for the Ricci flow (Cf. [H4])
since it implies a Harnack type estimate.
Since this  and similar results originate from the fundamental work
of Li-Yau \cite{L-Y} and Hamilton \cite{H4}, it seems to be
more appropriate to call them Li-Yau-Hamilton type inequalities. We shall
adopt this terminology in this work.

The second author would like to thank Shing-Tung Yau for useful conversations.

%Feb 28, 2001
\input amstex
\documentstyle{amsppt}
\magnification=1200
\hsize=13.8cm
\def \b {\beta}
\def\i{\sqrt{-1}}
\def\Ric{\text{Ric}}
\def\lf{\left}
\def\ri{\right}
\def\bbar{\bar \beta}
\def\a{\alpha}
\def\g{\gamma}
\def\p{\partial}
\def\delbar{\bar\delta}
\def\ddbar{\partial\bar\partial}

\def\C{\Bbb C}

\def\vp{\varphi}
\def\tD{\tilde \Delta}

\def\bb{{\bar\beta}}
\def\abb{{\alpha\bar\beta}}

\def\i{\sqrt {-1}}
\def\tD{\widetilde \Delta}
\def\tn{\widetilde \nabla}
\def \D {\Delta}

\def\aint{\frac{\ \ }{\ \ }{\hskip -0.4cm}\int}
\subheading{\S 1 Long time existence via Poincar\'e-Lelong equation}

Let $(M^m, g_{\a\bbar}(x))$ be a complete noncompact
K\"ahler manifold with bounded nonnegative holomorphic bisectional curvature.
Consider the K\"ahler-Ricci flow:
$$
\frac{\p}{\p t}g_{\a\bbar}(x,t)=-R_{\a\bbar}(x,t)\tag 1.1
$$
such that $g_{\a\bbar}(x, 0)=g_{\a\bbar}(x)$.

In \cite{Sh1-3}, short time existence of (1.1) was established, and  the long time existence was also proved under the assumption that
$$
\aint_{B(x,r) }\Cal R_0 dV\le Cr^{-\theta}\tag1.2
$$
for some constants $C$ and $\theta>0$ for all $x$ and  $r$. Here $\Cal R_0$ is the scalar curvature of the initial metric and $\aint_{B_x(r)}\Cal R_0 dV$ is the average of $\Cal R_0$ on the geodesic ball $B(x,r)$ with center at $x$ and radius $r$. The proof of the long time existence in \cite{Sh2, Sh3} is rather complicated. In this section, with the help of solutions of the Poincar\'e-Lelong equation we shall give a simple proof of the long time existence by using a maximum principle. Our assumption on $\Cal R_0$ is a little bit different from
(1.2).

Let us recall the result on short time existence of
Shi \cite{Sh3}.

\proclaim{Theorem 1.1} Let $(M^m, g_{\a\bbar}(x))$ be a complete noncompact K\"ahler
manifold with   nonnegative holomorphic bisectional curvature such that the
scalar curvature $\Cal R_0$ is bounded by $C_0$. Then (1.1) has a solution on $M\times[0,T)$ for some $T>0$ depending only $m$ and $C_0$ such that the following are true.
\roster
\item"{(i)}" $(M,g_\abb(x,t))$ is a K\"ahler metric with nonnegative holomorphic bisectional curvature for $0\le t<T$.
\item"{(ii)}" There exists $C>0$ such that
$$
C^{-1}g_\abb(x,0)\le g_\abb(x,t)\le g_\abb(x,0),\tag1.3
$$
and
$$
0\le\Cal R(x,t)\le C\tag1.4
$$
for all $(x,t)\in M\times[0,T)$.
\endroster
\endproclaim

Before we give our proof on the long time existence, let us fix the
 notations.  For any smooth function $f$, let
$\Delta f=g^{\a\bbar}(x,t)\frac{\p ^2 f}{\p z^\a \p \bar{z}^{\beta}}$, $|\nabla f|^2=g^{\a\bbar}(x,t)f_{\a}f_{\bbar}$. Summation  convention is understood. We also use $\tD$
and $\tn$ to denote the Laplacian and the gradient with respect to
a fixed metric  $g_{\a\bbar}(x)$ or the initial metric $g_{\a\bbar}(x,0)$ of the solution of (1.1).  $B_t(x,r)$ is the geodesic ball of radius $r$ with  respect to the metric $g_{\a\bbar}(x,t)$ and $V_t(x,r)$ be the volume of
$B_t(x,r)$ with respect to $g_{\a\bbar}(x,t)$.
We may also
use the ones without $t$  to denote the balls and volumes for a fixed metric. The same convention  applies to the  distance function $r_t(x,y)$ between two points $x,y\in M$ as well as the  volume element $dV_t$. As in [Sh2], throughout this work, let
$$
F(x,t)=\log\left(\frac{\det(g_{\a\bbar}(x,t))}
{\det(g_{\a\bbar}(x,0))}\right).
$$
Then for the solution of (1.1)
$$
dV_t=e^FdV,\tag1.5
$$
$$
F(x,t)=-\int_0^t\Cal R(x,\tau)d\tau\tag1.6
$$
where $\Cal R(x,t)$ is the scalar curvature of the metric $g_\abb(x,t)$. For the solution of  (1.1), we have the following maximum principle, which is of independent interest. The proof follows the idea in \cite{K-L}
(see also Li's lecture  notes [Li]).

 Let  $g_{ij}(x,t)$ be a smooth family of complete Riemannian metrics defined on $M$ with $0\le t\le T_1$ for some $T_1>0$ satisfying  the following properties: There exists a constant $C_1>0$ such that for any $T_1\ge t_2\ge t_1\ge0$
$$
C_1 g_{ij}(x,t_1)\le g_{ij}(x,t_2)\le g_{ij}(x,t_1)\tag1.7
$$
for all $x\in M$.

\proclaim{Theorem 1.2} With the above assumptions and notations, let $f(x,t)$ be a smooth function such that
$(\Delta -\frac{\p}{\p t})f(x,t)\ge 0$ whenever $f(x,t)\ge0$. Assume that
  $$
\int_0^{T_1}\int_M \exp(-ar_0^2(x))f_{+}^2(x,s)\, dV_0\, ds <\infty \tag 1.8
$$
for some $a>0$, where $r_0(x)$ is the distance function
to a fixed point $o\in M$ with respect to  $g_{ij}(x,0)$. Suppose $f(x,0)\le 0$ for all $x\in M$. Then $f(x,t)\le 0$ for all $(x,t)\in M\times [0,T_1]$.
\endproclaim
\demo{Proof}  Let $F(x,t)$ be such that $dV_t=e^F(x,t)dV_0$. By (1.7), we have
$$
\frac{\p}{\p t}F\le 0.\tag1.9
$$
Let $0<T\le T_1$ which will   be specified later and
 let
$$
g(x,t)=\frac{-r^2_{T}(x)}{4(2T-t)}, \ \ \ \text{  on } M\times [0,T].
$$
Here $r_{T}(x)$ is the distance function to $o\in M$ with respect to
$g_{\a\bbar}(x,T)$.  It is  easy to
check that
$$
|\nabla_{T} g|^2+\frac{\p g}{\p t} =0
$$
Here $\nabla_{T}$ is the gradient with respect to $g_{\a\bbar}(x,T)$. By (1.7), $g_{ij}$ is nonincreasing in $t$, hence we have
$$
|\nabla g|^2+\frac{\p g}{\p t} \le |\nabla_{T} g|^2+\frac{\p g}{\p t}= 0, \tag 1.10
$$
for $t\in [0,T]$. Let $\varphi(x)$ be a cut-off function  which we will specify later.  We have
$$
\split
0& \le \int_0^{T} \int_M \vp^2 e^g f_+\left(\Delta-\frac{\p}{\p t}\right)f
\, dV_s\, ds \\
& =\int_0^{T} \int_M \vp^2 e^g f_+ (\Delta f) \, dV_s\, ds -\frac{1}{2}
\int_0^{T} \int_M \vp^2 e^g \frac{\p}{\p t}(f_+ ^2)\, dV_s\, ds.
\endsplit
\tag1.11
$$
Here $f_+:=\max\{0, f\}$.
Now we calculate the last two terms in the above inequality.
$$
\split
\int_M \vp^2 e^g f_+ (\Delta f) \, dV_s & = -\int_M\vp^2e^g|\nabla f_+|^2 \,
dV_s
-2\int_M \vp e^g <\nabla \vp, \nabla f_+>f_+\, dV_s\\
&\ \ \ -
\int_M \vp^2e^g f_+<\nabla g, \nabla f_+>\, dV_s \\
&\le 2\int_M e^g f_+ ^2 |\nabla \vp|^2\, dV_s+\frac{1}{2}\int_M \vp^2e^gf_+ ^2|\nabla g|^2\, dV_s.
\endsplit
\tag1.12
$$
On the other hand,
$$
\split
-\frac{1}{2}\int_0^{T}\int_M \vp^2 e^g\frac{\p}{\p t}(f_+ ^2)\, dV_s\, ds
& = \frac{1}{2}\bigg[-\int_M \vp^2 e^g f_+ ^2\, dV_s\bigg|^{T}_{0} +
\int_0^{T}\int_M \vp^2 e^g g_s f_+ ^2\, dV_s\, ds \\
&\ \  +\int_0^{T}\int_M \vp^2e^gf_+ ^2F_s(y,s)\, dV_s\, ds\bigg]\\
&\le \frac{1}{2}\bigg[-\int_M \vp^2 e^g f_+ ^2\, dV_s\bigg|^{T}_{0} +
\int_0^{T}\int_M \vp^2 e^g g_s f_+ ^2\, dV_s\, ds \bigg]
\endsplit
\tag 1.13
$$
where we have used (1.9). Combining (1.10)--(1.13),  we have that
$$
\int_M \vp^2(x)e^{g(x,T)}f_+^2(x,T)\, dV_T  \le 4\int_0^{T}\int_M e^gf_+ ^2|\nabla \vp|^2\,
dV_s\, ds.
$$
Now using (1.7)  we have
$$
\int_M \vp^2(x)e^{g(x,T)}f_+^2(x,T)\, dV_T \le C_3\int_0^{T}\int_M
e^gf_+^2|\tn\vp|^2\, dV_0\, ds  \tag 1.14
$$
for some constant $C_3$ depending on $C_1$ in (1.7).
Here $\tn$ is the gradient with respect the initial metric
$g_{ij}(x,0)$.
For $R>0$, let $\vp$ be the function with compact support such that
$$
\split
&\vp(x)=1,\ \ \ \text{for } \ x\in B_0(o,R);\\
&\vp(x)=0, \ \ \ \text{for }\ x\in M\setminus B_0(o,2R);\\
&|\tn \vp|\le \frac{2}{R}.
\endsplit
$$
Letting $R\to \infty$ in (1.14) we have that
$$
\int_M  e^{g(x,T)}f_+^2(x,T)\, dV_T \le \liminf_{R\to \infty}\frac{4C_3}{R^2}
\int_0^{T}
\int_{B_0(o,2R)\setminus B_0(o,R)}e^{-\frac{r^2_0(x)}{C_4T}} f_+ ^2\, dV_0\, ds
$$
for some constant $C_4>0$ depending only on $C_1$ in (1.7).
Now if  $T<\frac{1}{aC_4}$, by (1.9), we will have  $$
\int_M  e^{g(x,T)}f_+^2(x,T)\, dV_T \le 0.
$$
This implies that $f(x,T)\le 0$. Since $C_4$ depends only on $C_1$,
iterating this procedure we complete the proof of the theorem.
\enddemo

Let $g_\abb(x,t)$ be a solution of (1.1) on $M\times [0,T)$, which is K\"ahler for all $t$. We have the following easy lemma.
\proclaim{Lemma 1.1} Suppose there is a function $u_0(x)$ such that
$$
\i\ddbar u_0 =\Ric(g(\cdot,0))\tag1.15
$$
where $\Ric(g(0))$ is the Ricci form of the initial metric $g(0)$. Let $F$ be the ratio of the volume element as in (1.5) and let $u(x,t)=u_0(x)-F(x,t)$. Then
$$
\i\ddbar u =\Ric(g(t)),\tag1.16
$$
$$
\lf(\Delta-\frac{\p}{\p t}\ri)u(x,t)=0,\tag1.17
$$
 $$
\left(\Delta-\frac{\p}{\p t}\right)|\nabla u|^2 = \|u_{\a\beta}\|^2+
\|u_{\a\bbar}\|^2, \tag 1.18
$$

$$
\left(\Delta-\frac{\p}{\p t}\right)\lf(|\nabla u|^2+1\ri)^\frac12\ge0,\tag1.19
$$
and
$$
\left(\Delta-\frac{\p}{\p t}\right)\Cal R=\left(\Delta-\frac{\p}{\p t}\right)u_t=-\|u_{\a\bbar}\|^2.\tag 1.20
$$
Here $\|u_{\a\bbar}\|^2(x,t)=g^{\a \bbar}(x,t)g^{\g\delbar}(x,t)
u_{\a \delbar}(x,t)u_{\g\bbar}(x,t)$,  $ \|u_{\a\b}\|^2(x,t)$
$=g^{\abb}(x,t)$
$g^{\g\delbar}(x,t)$ $ u_{\a\gamma}(x,t)$ $ u_{\bbar\delbar}(x,t)$.
\endproclaim
\demo{Proof} (1.16) and (1.17)   follow from the fact that $g_\abb(x,t)$ is a solution of (1.1) which is K\"ahler, and the definition of $F$ and $u_0$.

To prove (1.18), after choosing a normal  coordinates with respect to $g_\abb(x,t)$ near any fixed point
$$
\split
\D \, |\nabla u|^2 &  =g^{\g\bar{\delta}}
\left(u_\a u_{\bar{\b}}g^{\a\bar{\b}}\right)
_{\g\bar{\delta}}\\
& =  u_{\a\g}u_{\bar{\a}\bar{\g}}+u_{\a\bar{\g}}u_{\bar{\a}\g}
+
(\D u)_\a u_{\bar{\a}}+u_\a (\D u)_{\bar{\a}}+u_{\a\bar{\b}}u_\a u_{\bar{\b}},
\endsplit
$$
where we have used (1.1) and (1.16). Using (1.1), we have
$$
\frac{\p}{\p t}\, |\nabla u|^2\,  =\,  (u_t)_{\a}u_{\bar{\a}} + u_\a(u_t)_{\bar{\a}}
+u_{\a\bar{\b}}u_a u_{\bar{\b}}.
$$
Combining this with (1.17), we have (1.18). (1.19) follows from (1.18) by direct computations.

To prove (1.20), differentiate (1.17) with respect to $t$. Using (1.16) we have
$$
\split
\lf(\D-\frac{\p}{\p t}\ri)\Cal R &=\lf(\D-\frac{\p}{\p t}\ri)u_t\\
&=-g_t^{\abb}u_\abb\\
&=g^{\xi\bb}g^{\a\bar\g}g_{\xi\bar\g,t}u_\abb\\
&=-g^{\xi\bb}g^{\a\bar\g}R_{\xi\bar\g}u_\abb \\
&=-g^{\xi\bb}g^{\a\bar\g}u_{\xi\bar\g}u_\abb\\
&=-\|u_\abb\|^2.
\endsplit
$$
This completes the proof of the lemma.
\enddemo

We are ready to prove the long time existence.
\proclaim{Theorem 1.3} Let $(M^m, g_\abb(x,t))$ be a complete noncompact K\"ahler manifold with nonnegative holomorphic bisectional curvature such that its scalar curvature $\Cal R_0$ is bounded and
satisfies
$$
\int_0^\infty k(x,s)ds\le C_1\tag1.21
$$
for some constant $C_1$ for all $x$ and $r$, where
$$
k(x,s)=\aint_{B(x,s)}\Cal R_0dV.
$$
Then (1.1) has long time existence. Moreover, there is a function $u(x,t)$ such that
$$\i\ddbar u(\cdot,t)=\Ric(g(t)),\tag1.22
$$
$$
|\nabla u|\le C(m)C_1,\tag1.23
$$
and
$$
\Cal R(x,t)+|\nabla u|^2(x,t)\le
\sup_{x\in M}\lf(\Cal R_0(x)+|\tn u_0|^2(x)\ri)\le \sup_{x\in M}\Cal R_0(x)+\lf(C(m)C_1\ri)^2\tag1.24
$$
for some constant positive $C(m)$ depending only on $m$ and for all $(x,t)$.
Moreover, the equality  holds for some $(x_0,t_0)$, with $t_0>0$ if and only if
$g_{\abb}(x,t)$ is a K\"ahler-Ricci soliton.
\endproclaim
\demo{Proof}
 By Theorem 1.1, there is a maximal $\infty\ge T_{\text{\rm max}}>0$ such that (1.1) has a solution $g_\abb(x,t)$ which satisfies condition (i) in Theorem 1.1 for   $0\le t<T_{\text{\rm max}}$, and satisfies the following condition: For any $0<T<T_{\text{\rm max}}$, there is a constant $C>0$ such that (1.3) and (1.4) are true  on $M\times [0,T]$.  By (1.21) and the results in \cite{N-S-T, Theorems 1.3 and 5.1}, there is a function $u_0(x)$ such that
$$
\i\ddbar u_0=\Ric (g(0))
$$
and
$$
|\tn u_0|(x)\le C(m)C_1\tag1.25
$$
for all $x$ for some constant $C(m)$ depending only on $m$. Let  $u(x,t)=u_0(x)-F(x,t)$ and let $0<T<T_{\text{\rm max}}$ be fixed. By (1.3), (1.4), (1.6) and (1.25), it is easy to see that there is a constant $C_2$ such that for $(x,t)\in M\times[0,T]$
$$
|u(x,t)|\le C_2(r_0(x)+1)\tag1.26
$$
where $r_0(x)$ is the distance from a fixed point $o$ with respect to  $g(0)$. By Lemma 1.1 (1.16), we have  $\Delta u(x,t)=\Cal R(x,t)$. Combining this
with   (1.4) and (1.26), it is not hard to prove that
$$
\int_{B_t(o,r)}|\nabla u|^2\le C_3r^{2m+1}\tag1.27
$$
for some constant $C_3$ for all $0\le t\le T$ and for all $r$. Here we have
 used the fact that $g_\abb(x,t)$ has nonnegative Ricci curvature and volume
 comparison. Hence using (1.3), we conclude that the function
 $f=\lf(|\nabla u|^2+1\ri)^\frac12-\lf(C^2(m)C^2_1+1\ri)^\frac12$ satisfies
 the condition (1.10) in
Theorem 1.2 with $T_1$ replaced by $T$. Here $C(m)$ is the constant in (1.25).
 By (1.19) of Lemma 1.1 and Theorem 1.2, we can conclude that (1.23) is true
 for $x\in M$ and  $0\le t\le T_{\text{\rm max}}$, because $T$ can be any
positive number less than $T_{\text{\rm max}}$.

By (1.18) and (1.20) of Lemma 1.1, we have
$$
\lf(\D-\frac{\p}{\p t}\ri)\lf(|\nabla u|^2+\Cal R\ri)= \|u_{\a\b}\|^2.\tag1.28
$$
By (1.23) and (1.4), we conclude that $|\nabla u|^2+\Cal R$ is uniformly bounded on $M\times[0,T]$. By (1.28), we can apply Theorem 1.2 again and conclude that (1.24) is true for all
$x\in M$ and  $0\le t\le T_{\text{\rm max}}$.
In particular $\Cal R$ is uniformly bounded on $M\times [0,  T_{\text{\rm max}})$. By Theorem 1.1, $T_{\text{\rm max}}$ must be infinity.
If for some $(x_0, t_0), t_0>0$,
$$
\left({\Cal R}+|\nabla u|^2\right)(x_0, t_0) =
\sup_{x\in M}({\Cal R}+|\nabla u|^2)(x,0)
$$
we can conclude that ${\Cal R}(x,t)+|\nabla u|^2(x,t)$ is constant,
by the strong maximum
principle. Thus $u_{\a\b}(x,t)=0$ by (1.28). Together with the fact
$u_{\abb}(x,t)=R_{\abb}(x,t)$, it implies that $g_{\abb}(x,t)$ is a K\"ahler-Ricci
soliton. It is easy to check that for a K\"ahler Ricci soliton (1.24) holds with
the equality (Cf. [C-H]).

\enddemo
%\end

%Feb 28, 2002
\input amstex
\documentstyle{amsppt}
\magnification=1200
\hsize=13.8cm
\def \b {\beta}
\def\i{\sqrt{-1}}
\def\Ric{\text{Ric}}
\def\lf{\left}
\def\ri{\right}
\def\bbar{\bar \beta}
\def\a{\alpha}
\def\g{\gamma}
\def\p{\partial}
\def\delbar{\bar\delta}
\def\ddbar{\partial\bar\partial}

\def\C{\Bbb C}

\def\vp{\varphi}
\def\tD{\tilde \Delta}

\def\bb{{\bar\beta}}
\def\abb{{\alpha\bar\beta}}

\def\i{\sqrt {-1}}
\def\tD{\widetilde \Delta}
\def\tn{\widetilde \nabla}
\def \D {\Delta}
\def\R{\Bbb R}

\def\aint{\frac{\ \ }{\ \ }{\hskip -0.4cm}\int}
\subheading{\S 2 Some properties preserved by the K\"ahler-Ricci flow}

In this section, we shall investigate the behavior of
 $\aint_{B_t(x_0,r)}\Cal RdV_t$. To do this, we shall give some generalizations
 of the estimates in \cite{Sh2-3, C-Z, C-T-Z} from above and below on
the volume element $F(x,t)$ defined in (1.5). More precisely,
we shall obtain upper and lower estimates on $F(x,t)$ in terms of the integral
$$
\int_0^rsk(x,s)ds
$$
where $k(x,s)$ is the average of the scalar curvature $\Cal R_0$ over $B_0(x,s)$ at $t=0$. Our proofs use the well-known estimates of the heat kernels
and the Green's functions for manifolds with nonnegative Ricci curvature
of Li-Yau \cite{L-Y}.  Our proofs seem to be simpler than those in
[Sh2-3], etc. Also we do not use the complicated construction of exhaustion
 functions as in the \cite{Sh2-3, C-Z, C-T-Z}.
To derive our estimates we need the following lemma,
which is a direct consequence of the mean value inequality
of Li-Schoen \cite{L-S} on subharmonic functions.

\proclaim{Lemma 2.1 (Generalized mean value inequality)}  Let $M^n$ be a complete noncompact Riemannian manifold with nonnegative Ricci curvature with real dimension $n$.  Let $u\ge 0$ be a smooth function such that  $ \tD u\ge -f$ with $f\ge0$. For any $x_0\in M$ and $r>0$, we have
$$
u(x_0)\le \int_{B(x_0,r)}G_r(x_0,y)f(y)dy+C(n)\aint_{B(x_0,r)}u\tag2.1
$$
for some constant $C(n)$ depending only on $n$, where $G_r(x,y)$ is the positive Green's function on $B(x_0,r)$ with zero boundary value.
\endproclaim
\demo{Proof} Let $v$ be such that $\tD v=-f$ on $B(x_0,r)$ and $v=0$ on $\partial B(x_0,r)$. Note that $v\ge0$ in $B(x_0,r)$.  Since $w=\max\{u-v,0\}$ is Lipschitz, subharmonic and nonnegative,   by the mean value inequality of Li-Schoen \cite{L-S}, we have
$$
w(x_0)\le C  \aint_{B(x_0,r)}w
$$
for some constant $C=C(n)$ depending only on $n$. If $u(x_0)-v(x_0)\le 0$, then we have
$$
u(x_0)\le v(x_0)=\int_{B(x_0,r)}G_r(x_0,y)f(y)dy.
$$
In this case, (2.1) is true. If $u(x_0)-v(x_0)>0$ then
$$
\split
u(x_0)
&=w(x_0)+v(x_0)\\
&\le C\aint_{B(x_0,r)}w+v(x_0)\\
&\le C\aint_{B(x_0,r)}u+v(x_0)\\
&\le C\aint_{B(x_0,r)}u+\int_{B(x_0,r)}G_r(x_0,y)f(y)dy
\endsplit.
$$
Therefore (2.1) is also true for this case.
\enddemo

We should mention that the above lemma was also proved in a somewhat
different form in [Sh2-3] with a more complicated  proof (Cf.
Lemma 6.10 of [Sh2] and Lemma 6.8 of [Sh3]). We also need
the following estimates of Green's functions.
\proclaim{Lemma 2.2} Let $M^n$ be as in Lemma 1.1. For any function
 $f\ge0$,   let $k(x,r)=\aint_{B(x,r)}f$. Then we have
$$
\int_{B(x,r)}G_r(x,y)f(y)dy\ge C(n)\lf(r^2k(x,\frac r5)+\int_0^{\frac r5}sk(x,s)dr\ri),
$$
for some constant $C(n)>0$ depending only on $n$, where $G_r$ is the Green's function on $B(x,r)$ where zero boundary value. If in addition, $M$ supports a minimal positive Green's function $G(x,y)$ such that
$$
\alpha\cdot\frac{r^2(x,y)}{V(x, r(x,y))}\le G(x,y)\le \frac1\alpha\cdot \frac{r^2(x,y)}{V(x, r(x,y)}.
$$
for some $\alpha>0$ for all $x,y\in M$, then

$$
\int_{B(x,r)}G(x,y)f(y)dy\le C(n,\alpha)\lf(r^2k(x,r)+\int_0^rsk(x,s)dr\ri),
$$
for some positive constant $C(n,\alpha)$ depending only on $n$ and $\alpha$.
\endproclaim
\demo{Proof} See the proofs of \cite{N-S-T, Theorems 1.1, 2.1}.
\enddemo

In the rest of this section, we assume $M^m$ is a complete noncompact K\"ahler manifold with bounded nonnegative holomorphic bisectional curvature such that $g_\abb$ is a solution of (1.1) on $M\times[0,T)$ with $T\le \infty$. We also assume that conditions (i) and (ii) are satisfied by $g_\abb$ on $M\times[0,T_1]$ for any $T_1<T$. Let $\frak m(t)=\inf_MF(\cdot,t)$. Then $\frak m(t)\le 0$.

With the notations as in \S1,  we also  need the following result
of Shi \cite{Sh3, p. 156}.
\proclaim{Lemma 2.3}
$$
\split
\Cal R_0(x)&\ge\Cal R_0(x)+e^FF_t\\
&\ge \Cal R_0(x)-g^{\alpha\bar\beta}(x,0)R_{\alpha\bar\beta}(x,t)\\
&=\tD  F(x,t)\\
&\ge \Cal R_0(x)-\Cal R(x,t)
\endsplit
\tag2.2
$$
where $\tD  $ is the Laplacian of the metric $g(0)$.
\endproclaim
\proclaim{Theorem 2.1} With the above assumptions and notations, the
 following estimates are true. Namely there exists $C_1>0$ depending only on
 $m$ such that for all $(x_0,t)\in M\times[0,T)$
$$
-F(x_0,t)\ge C_1^{-1}\int_0^{\sqrt t}sk(x_0,s)ds\tag2.3
$$
and
$$
-F(x_0,t)\le C_1\lf[\lf(1+\frac{t\lf(1-\frak m(t)\ri)}{R^2}\ri)\int_0^{R}sk(x_0,s)ds-\frac{t\frak m(t)\lf(1-\frak m(t)\ri)}{R^2}\ri],
\tag2.4$$
where $k(x_0,t)=\aint_{B_0(x_0,r)}\Cal R_0dV_0$.
\endproclaim
\demo{Proof} To prove (2.3), by Lemma 2.3 we have
$$
\tD F\ge \Cal R_0-\Cal R=\Cal R_0+F_t
$$
and so
$$
\lf(\tD-\frac{\p}{\p t}\ri)(-F)\le -\Cal R_0.\tag2.5
$$
Let $H(x,y,t)$ be the heat kernel of $M$ with respect to the metric $g(0)$, and let
$$
v(x,t)=\int_0^t\int_MH(x,y,t)\Cal R_0(y)dV_0(y).
$$
Then $\tD v-v_t=-\Cal R_0$ and $v=0$ at $t=0$. By (2.5) and the fact that $F(\cdot,0)\equiv0$, by the maximum principle and the estimate of the heat kernel \cite{L-Y}, we have for $(x,t)\in M\times[0,T)$
$$
\split
-F(x,t) &\ge v(x,t)\\
&=\int_0^t\int_M H(x,y,\tau){\Cal R}_0(y)\, dV_0\, d\tau \\
       & \ge C_2\int_0^t\int_0^\infty \frac{1}{V_0(x,\sqrt{\tau})}
e^{-\frac{r^2}{5\tau}}\int_{\partial B_0(x,r)}{\Cal R}_0(y)\, dA_0 \, dr\, d\tau\\
&\ge C_2\int_0^t\int_0^{\sqrt{\tau}}\frac{1}{V_0(x,\sqrt{\tau})}
e^{-\frac{r^2}{5\tau}}\int_{\partial B_0(x,r)}{\Cal R}_0(y)\, dA_0 \, dr\, d\tau\\
&= C_3 \int_0^t k(x,\sqrt{\tau})\, d\tau\\
& = 2C_3\int_0^{\sqrt{t}} \tau k(x,\tau)\, d\tau.
\endsplit
$$
for some positive constants $C_2-C_3$ depending only on $m$. Hence (2.3) is true.

To prove (2.4), by Lemma 2.3, $\tD F\le \Cal R_0+e^FF_t$. Hence for any  $(x_0,t)\in M\times[0,T)$ for any $R>0$, integrating the above inequality over $B_0(x_0,R)\times[0,t]$, we have
$$
\split
\int_0^t\int_{B_0(x_0,R)}&G_R(x_0,y)\tD F(y,s)dV_0ds\\
&\le t\int_{B_0(x_0,R)}G_R(x_0,y) \Cal R_0(y)dV_0+\int_{B_0(x_0,R)}
G_R(x_0,y)(e^{F(y,t)}-1)  dV_0,
\endsplit
$$
and
$$
\split
\int_{B_0(x_0,R)}
&G_R(x_0,y)(1-e^{F(y,t)} )  dV_0\\
&\le t\int_{B_0(x_0,R)}G_R(x_0,y) \Cal R_0(y)dV_0+\int_0^t\int_{B_0(x_0,R)}G_R(x_0,y)\tD\lf(- F(y,s)\ri)dV_0.
\endsplit\tag 2.6
$$
By the Green's formula, for each $0\le s\le t$
$$
\split
\int_{B_0(x_0,R)}G_{R}(x_0,y)  \tD(-F(y,s))dV_0 &=
F(x_0,s)+
\int_{\partial B_0(x_0,R)}F(y,s)\frac{\partial G_{R}(x_0,y)}{\partial \nu} \\
&\le  -\frak m(t) \endsplit,
$$
where we have used  the fact that $\frak m(t)$ is nonincreasing, $F\le 0$,  $\frac{\partial }{\partial \nu}G_{R}(x_0,y)\le0$ and $\int_{\partial B_0(x_0,R)}\frac{\partial }{\partial \nu}G_{R}(x_0,y)=-1$. Combining this with (2.6), we have
$$
\int_{B_0(x_0,R)}
G_R(x_0,y) (1-e^{F(y,t)} )dV_0\le t\lf(\int_{B_0(x_0,R)}G_R(x_0,y) \Cal R_0(y)dV_0-\frak m(t)\ri).
$$
Using the first inequality in Lemma 2.2,  this implies
$$
R^2\aint_{B_0(x_0,\frac15R)}\lf (1-e^{F(y,t)} \ri)dV_0\le C_4t\lf(\int_{B_0(x_0,R)}G_R(x_0,y) \Cal R_0(y)dV_0-\frak m(t)\ri)
 \tag2.7
$$
for some constant $C_4$ depending only on $m$.
Since if $0\le x\le 1 $, $1-e^{-x}\ge \frac13x$, we have $(1-e^F)\lf(1-\frak m(t)\ri)\ge -CF$ for some absolute positive constant $C$.
 Hence (2.7) implies that
$$
R^2\aint_{B_0(x_0,\frac15R)}\lf (-F(y,t) \ri)dV_0\le  C_5t\lf(1-\frak m(t)\ri)\lf(\int_{B_0(x_0,R)}G_R(x_0,y) \Cal R_0(y)dV_0-\frak m(t)\ri)
 \tag2.8
$$
for some  constant $C_5$ depending only on $m$. By Lemma 2.3, $\tD(-F)\ge -\Cal R_0$. By Lemma 2.1 and (2.8), there is a constant $C_6$ depending only on $m$ such that
$$
\split
-F(x_0,t)&\le \int_{B_0(x_0,\frac15R)}G_{\frac15R}(x_0,y) \Cal R_0(y)dV_0+C(n)\aint_{B_0(x_0,\frac15R)}\lf(-F(y,t)\ri)dV_0\\
&\le \int_{B_0(x_0,\frac15R)}G_{\frac15 R}(x_0,y)\Cal R_0(y)dV_0\\
&\qquad+\frac{ C_6 t\lf(1-\frak m(t)\ri)}{R^2}\lf(\int_{B_0(x_0,R)}G_R(x_0,y) \Cal R_0(y)dV_0-\frak m(t)\ri),
\endsplit \tag2.9
$$
where $G_{\frac15R}$ is the Green's function on $B_0(x_0,\frac15R)$. As in \cite{Sh3}, by considering $M\times\C^2$, we may assume that $M$ has positive Green's function which satisfies the condition in Lemma 2.2. Applying Lemma 2.2, we can conclude from (2.9) that
$$
-F(x_0,t)\le C_7\lf[\lf(1+\frac{t\lf(1-\frak m(t)\ri)}{R^2}\ri)\int_0^{2R}sk(x_0,s)ds-\frac{t\frak m(t)\lf(1-\frak m(t)\ri)}{R^2}\ri],
$$
for some constant $C_7$ depending only on $m$. This completes the proof of the theorem.
\enddemo
\proclaim{Corollary 2.1} Same assumptions and notations as in Theorem 2.1. Suppose   $k(x,r)\le k(r)$ for some function $k(r)$ for all $x\in M$. Then there exist positive constants $C$, $a$ depending only on $m$ such that for $0\le t<T$
$$
-\frak m(t)\le C\int_0^Rsk(s)ds\tag2.10
$$
where $R^2= at(1-\frak m(t))$.
\endproclaim
\demo{Proof} By (2.4), we have for any $R>0$
$$
-\frak m(t)\le C_1 \lf[\lf(1+\frac{t\lf(1-\frak m(t)\ri)}{R^2}\ri)\int_0^{ R}sk(s)ds-\frac{t\frak m(t)\lf(1-\frak m(t)\ri)}{R^2}\ri]
$$
where $C_1$ is a constant depending only on $m$. Let $R^2=2C_1t(1-\frak m(t))$, we have
$$
-\frak m(t)\le 2C_1\lf(1+\frac1{2C_1}\ri)\int_0^Rsk(s)ds.
$$
From this the result follows.
\enddemo
\proclaim{Corollary 2.2} With the same assumptions as in Corollary 2.1. Suppose
$$
\int_0^rsk(s)ds\le r^2\phi(r)
$$
for all $r$, where $\phi(r)$ is a nonincreasing function of $r$ such that $\lim_{r\to\infty}\phi(r)=0$. For $0<\tau\le \sup\phi$, let
$$
\psi(\tau)=\sup\{r|\ \phi(r)\ge \tau\}.
$$
Then for $0\le t<T$,
$$
-\frak m(t)\le \max\{1,\frac {C'}t\psi^2(\frac {C''}t)\}
$$
for some positive constants $C'$ and $C''$ depending only on $m$. In particular, the  K\"ahler-Ricci flow has long time existence.
\endproclaim
\demo{Proof} Note the $\psi(\tau)$ is finite  and nonincreasing for $0<\tau\le \sup\phi$ because $\phi(r)\to0$ as $r\to\infty$. By Corollary 2.1, there exist constants $a$ and $C_1$ depending only on $m$ such that
$$
-\frak m(t)\le C_1\int_0^Rsk(s)ds\le C_1at(1-\frak m(t))\phi\lf(\sqrt{at(1-\frak m(t))}\ri)
$$
where $R^2=at(1-\frak m(t))$. Suppose $-\frak m(t)\ge1$, then the above inequality implies that
$$
\phi\lf(\sqrt{at(1-\frak m(t))}\ri)\ge \frac{1}{2C_1at}.
$$
In particular, $\frac{1}{2C_1at}\le \sup \phi$. Hence
$$
\sqrt{at(1-\frak m(t))}\le \psi\lf(\frac{1}{2C_1at}\ri).
$$
Hence
$$
-\frak m(t)\le \max\{1,\frac {C'}t\psi^2(\frac {C''}t)\}
$$
for some positive constants $C'$ and $C''$ depending only on $m$.

The last statement follows from the method in \cite{Sh3, \S7}. Here we cannot use the method in Theorem 1.3 because we do not have a good solution for the Poincar\'e-Lelong equation.
\enddemo
\proclaim{Remark 2.1} The condition for long time existence in the corollary is weaker than that in \cite{Sh3}. In \cite{C-T-Z}, the long time existence is proved for the case of surfaces under the assumptions that the surface has maximal volume growth and that $\int_0^rsk(x_0,s)=o(r^2)$. The last assumption is a little bit weaker than ours.
\endproclaim
\proclaim{Remark 2.2} By the corollary, we may have the estimates in \cite{Sh2-3}. For example, if $k(r)=C(1+r)^{-2}$, then it is easy to see that $-\frak m(t)\le C\log (t+1)$. If $k(r)=C(1+r)^{-\theta}$ for $0<\theta<2$, then $-\frak m(t)\le C(t+1)^{(2-\theta)/\theta}$. In addition to these results in \cite{Sh2-3}, we may have the following estimate. Namely, if $\int_0^\infty k(r)dr<\infty$, then $-\frak m(t)=o(t)$ and  if $\int_0^rsk(s)ds\le Cr^2/\log(2+r)$, then we have $-\frak m(t)\le  e^{Ct}$ for some $C>0$.
\endproclaim

Another application of the corollary is a slight generalization of a gap theorem of Chen-Zhu \cite{C-Z}. In \cite{C-Z}, it is proved that if $M$ is a complete K\"ahler manifold with bounded nonnegative holomorphic bisectional curvature such that
$$
k(x_0,r)=\aint_{B_0(x,r)}\Cal R_0dV_0\le \epsilon (r)r^{-2}
$$
for all $x$ and $r$, where $\epsilon(r)\to0$ as $r\to\infty$. Then $M$ must be flat. Note that under this condition, the K\"ahler-Ricci flow has long time solution such that $R(x,t)$ is uniformly bounded on $M\times[0,\infty)$ by Theorem 1.3 and so $-\frak m(t)\le Ct$. Moreover
$$
\int_0^rsk(x,s)ds=o(\log r)
$$
uniformly.

 Using Corollary 2.1, we have:
\proclaim{Corollary 2.3} Let $(M^m, g)$ be complete K\"ahler manifold with bounded nonnegative holomorphic bisectional curvature such that the K\"ahler-Ricci flow (1.1) has long time solution.
\roster
\item"{(a)}" Suppose $M$ is nonflat and
 $-\frak m(t)\le Ct^k$ for some constant $C$ and $k>0$.  Then
$$
\liminf_{r\to\infty}\frac{\int_0^rsk(x,s)ds}{\log r}   >0,\tag2.11
$$
$$
\liminf_{t\to\infty}\frac{-F(x,t)}{\log t}>0,\tag2.12
$$
and
$$
\liminf_{t\to\infty}t\Cal R(x,t) >0,\tag2.13
$$
for all $x$, where $k(x,s)=\aint_{B_0(x,r)}\Cal R_0dV_0$.
\item"{(b)}" If the Poisson equation $\tD u=\Cal R_0$ has a  solution $u$
which is bounded from above, then $M$ is flat. In particular, any bounded from
above solution is a constant.
\endroster
\endproclaim
\demo{Proof} Note that if (2.11) is true for some $x$,
it is true for all $x$. Suppose $M$ is nonflat, then there exists $x_0$ such that $\Cal R_0(x_0)>0$. If  (2.11) is not true, then
there exists $R_i\to\infty$ such that
$$
 \int_0^{R_i}sk(x_0,s)ds\le \frac1i\log R_i.\tag2.14
$$
Let $t_i\to\infty$ be such that $  t_i (1-(\frak
m(t_i))^2 = R_i^2$. By (2.4), we have
$$
\split
-F(x_0,t_i)&\le  C_1\lf(\int_0^{R_i}sk(x_0,s)ds+1\ri)\\
&\le C_1\lf(\frac1i\log R_i+1\ri)\\
&\le C_2\lf(\frac1i\log t_i+1\ri)
\endsplit \tag2.15
$$
for some constants $C_1-C_2$   independent of $i$. Here we have used the assumption that $-\frak m(t)\le Ct^k$. We can then proceed as in \cite{C-Z}. For any $T>0$, by the Li-Yau-Hamilton type  inequality \cite{Co2-3}  for $t>T$,
$$
\frac{T}{t}\Cal R(x_0,T)\le \Cal R(x_0,t).
$$
Integrating from $T$ to $t_i$, we have
$$
T\log \frac{t_i}T\Cal R(x_0,T)\le -F(x_0,t_i)\le C_2\lf(\frac1i\log t_i+1\ri).
$$
Dividing both sides by $\log t_i$ and let $t_i\to\infty$, we have
$\Cal R(x_0,T)=0$. Since $T$ is arbitrary, we conclude that $\Cal R(x_0)=0$. This is a contradiction. Hence (2.11) is true.

If (2.12) is not true for some $x$, then by (2.3) in Theorem 2.1, (2.11) is not true for this $x$. Hence $M$ must be flat by the previous result.

By (2.12), for any $x\in M$ there exists $C_3>0$ and $t_0>0$ such that
$$
-F(x,t)\ge C_3\log t, \tag2.16
$$
for all $t\ge t_0$. By the Li-Yau-Hamilton type
 inequality in \cite{Co2-3}, for all $t>t_0$ and $s\le t$,
$$
\frac{t}{s}\Cal R(x,t)\ge \Cal R(x,s).
$$
Integrating over $s$ from $1$ to $t$ and using (2.16) we have
$$
\split
\lf(t\log t\ri)\Cal R(x,t)&\ge \int_1^t\Cal R(x,s)ds\\
&=-F(x,t)-\int_0^1\Cal R(x,s)ds\\
&\ge C_3\log t-\int_0^1\Cal R(x,s)ds.
\endsplit
$$
From this (2.13) follows.

 The proof of (b) follows from the proof of (a) and
Theorem 2.1 of [N-S-T].
\enddemo
\proclaim{Remark 2.3}  The argument above in fact also shows that {\it
any bounded solution to $\tD u={\Cal R}_0(x)$ is a constant} since if
$\tD u={\Cal R}_0(x)$ has a bounded solution, we then have long time solution
to (0.2) by Theorem 1.3 and Theorem 2.1 of [N-S-T].
In \cite{Cw},
a gradient estimate is obtained for the K\"ahler-Ricci flow under the assumption that there is a bounded potential function for the Ricci tensor.
If we assume the manifold has nonnegative holomorphic bisectional curvature, then this is only possible for flat manifolds.
\endproclaim

\proclaim{Corollary 2.4}
Same assumptions and notations as in Corollary 2.1. If we assume that
$$
\lim_{r\to \infty}\frac{\int_0^r sk(s)\, ds }{r}=0
$$
we have long time existence for the K\"ahler-Ricci flow with
$$\lim_{t\to \infty}\frac{-\frak m(t)}{t}=0$$
and
$$\lim_{t\to \infty} R(x,t)=0$$
uniformly for $x\in M$. If in addition,
we assume that
$(M, g(0))$ has maximum volume growth, $M$ is diffeomorphic to
$\R^{2m}$, in case $m\ge 3$ and homeomorphic to
$\R^{4}$, in case $m=2$. Moreover, $M$ is a Stein manifold.
\endproclaim
\demo{Proof}\ The first part just follows from Corollary 2.1 and the
Li-Yau-Hamilton type
inequality of Cao \cite{Co2-3} as in the proof of Corollary 2.3.
To prove that
$M$ is Stein and topologically $\R^{2m}$
one just need to use the observation that
the injectivity radius of $M$ has a uniform lower bound in the case of
the maximum volume growth and bounded curvature tensor. Also $|R(x,t)|\to
0 $, as $t\to \infty$, means that the K\"ahler-Ricci flow will
improves the injectivity radius to $\infty$ along the flow. The rest
argument is same as in section 3 of [C-Z].
\enddemo

Another corollary of the proof of Theorem 2.1  is a result on   the preservation of volume growth under the K\"ahler-Ricci flow. In \cite{Sh2} it was proved that the property of having maximum volume growth is preserved  under the assumption that ${\Cal R}_0(x)$ is of quadratic decay. In  \cite{C-Z, C-T-Z} it was generalized to the case of  more relaxed decay conditions on ${\Cal R}_0(x)$ using the same argument as [Sh2]. In \cite{H3}, it was proved  under the Ricci flow with nonnegative Ricci curvature, and under the stronger assumption that
the Riemannian curvature tensor of the initial metric goes to zero
pointwisely, then the volume ratio $\lim_{r\to\infty}r^{-n}V_t(r)$ is preserved.  In our case,  we have the following stronger result:
\proclaim{Theorem 2.2} With the same assumptions and notations as in Theorem 2.1. Suppose
$$
\int_0^rsk(x,s)ds=o(r^2)\quad \text{\rm  as $r\to\infty$}.
$$
Let $o\in M$ be a fixed point. Then for any $0<t<T$,
$$
\lim_{r\to\infty}\frac{V_t(o,r)}{V_0(o,r)}=1
$$
where $V_t(o,r)$ is the volume of the geodesic ball $B_t(o,r)$ with respect to the metric $g(t)$ for $0\le t<T$.
\endproclaim
\demo{Proof} Since $\Cal R(x,t)$ is uniformly bounded on $M\times[0,t]$, by Theorem 17.2 in \cite{H3},  $B_t(o, r)\subset B_0(o, r+C_1 t)$ for some constant $C_1$ independent of $r$.  Using the fact that $g(t)$ is nonincreasing in $t$,  we have that
$$
\split
V_t(o, r)  & \le V_t\lf(B_0(o, r+C_1 t)\ri)\\
&\le  V_0\lf(B_0(o, r+C_1 t)\ri)\\
&\le V_0(o,r)\cdot\lf(\frac{r+C_1t}{r}\ri)^{2m}.
\endsplit
$$
This implies that
$$
\limsup_{r\to\infty}\frac{V_t(o,r)}{V_0(o,r)}\le 1.
$$
Using the fact that $g(t)$ is nonincreasing in $t$ again, we have
$$
\split
V_t(o, r)& \ge \int_{B_0(o, r)}\, dV_t \\
 & = \int_{B_0(o, r)}e^{F(y,t)}\, dV_0\\
&=V_0(o, r)+\int_{B_0(o, r)}(e^{F(y,t)}-1)\, dV_0.
\endsplit \tag 2.17
$$
On the other hand, using (2.7) in the proof of Theorem 2.1 and using Lemma 2.2 as in the proof of (2.4),
we have
$$
 \aint_{B_0(o,r)}\lf (1-e^{F(y,t)} \ri)dV_0\le C_2r^{-2}t\lf(\int_0^{10r}sk(o,s)ds-\frak m(t)\ri)
$$
for some constant $C_2$ independent on $r$. Combining this with (2.17), we have
$$
\frac{V_t(o, r)}{V_0(o,r)}\ge 1- C_2r^{-2}t\lf(\int_0^{10r}sk(o,s)ds-\frak m(t)\ri).
$$
Since $\int_0^Rsk(s)ds=o(R^2)$, we have
$$
\liminf_{r\to\infty}\frac{V_t(o, r)}{V_0(o,r)}\ge 1.
$$
The theorem then follows.
\enddemo

It was proved in \cite{H3} that the condition $|Rm|\to0$ as $x\to\infty$ is preserved under the Ricci flow. Applying Theorem 2.1, we can prove the decay rate of the scalar curvature in the average sense is preserved under the K\"ahler-Ricci flow in a certain sense.
\proclaim{Theorem 2.3} Let $M^m$ be a complete noncompact K\"ahler manifold with bounded nonnegative holomorphic bisectional curvature. Suppose (1.1) has long time existence, such that for any $T>0$ the conditions (i) and (ii) in Theorem 1.1 are satisfied. Then the following are true:
\roster
\item"{(a)}" Suppose $\int_0^rsk(x,s)ds\le C(1+r)^{1-\epsilon}$ for some constants $C>0$ and $\epsilon>0$ for all $x$ and $r$. Then $\int_0^rsk_t(x,s)ds\le C'(1+r)^\delta$ where $\delta=\min\{1,2(1-\epsilon)/(1+\epsilon)\}$ for some constant $C'$ independent of $x, t, r$.
\item"{(b)}" Suppose $\int_0^rsk(x,s)ds\le C\log(r+2)$   for some constants $C>0$ for all $x$ and $r$. Then $\int_0^rsk_t(x,s)ds\le C'\log(r+2)$ for some constant $C'$ independent of $x,  t, r$.
\endroster
Here $k(x,r)=\aint_{B_0(x,r)}\Cal R_0dV_0$ and $k_t(x,r)=\aint_{B_t(x,r)}\Cal R(y,t)dV_t$.
\endproclaim
\demo{Proof} We  prove (b) first. For $T\ge0$,  let
$$
F(x,t;T)=\log\left[\frac{\det\lf(g_{\a\bbar}(x,t+T)\right)} {\det\left
(g_{\a\bbar}(x,T)\ri)}\right].
$$
Considering the flow $g_{\a\bbar}(x,t+T)$ with initial data   $g_{\a\bbar}(x,T)$ and using   (2.3) in Theorem 2.1, we have for any $t>0$
$$
-F(x,t;T)\ge C_1\int_0^{\sqrt t}sk_T(x,s)ds.\tag 2.18
$$
for some constant $C_1>0$ depending only on $m$. On the other hand,
by the Li-Yau-Hamilton  inequality \cite{Co2-3}
$$
T\Cal R(x,T)\le t\Cal R(x,t)
$$
for all $t\ge T$. We have
$$
\int_T^t\frac Ts\Cal R(x,T)ds\le \int_T^t\Cal R(x,s)ds\le -F(x,t;0)
\le C_2\log(t+2)
$$
for some constant $C_2$ independent of $x$ and $t$, where we have used Corollary 2.1 and the assumption on $k(x,r)$. Dividing both sides by $\log t$ and let $t\to\infty$, using the fact that $\Cal R$ is uniformly bounded on $M\times[0,\infty)$ by Theorem 1.3, we have
$$
 \Cal R(x,T) \le \frac{C_3}{T+1}\tag2.19
$$
for some constant $C_3$ independent of $x$ and $t$.   Since the metric is nonincreasing along the Ricci flow $\det\lf(g_{\a\bbar}(x,T)\ri)\le
 \det\lf(g_{\a\bbar}(x,0)\ri)$, by (2.18) and Theorem 2.1,    for all $t>0$
$$
\split
\log(t+T+2)&\ge -C_4F(x,t+T;0)\\
&\ge -C_4F(x,t;T)\\
&\ge C_5 \int_0^{\sqrt t}sk_T(x,s)ds
\endsplit\tag2.20
$$
for some positive constants $C_4-C_5$  independent of $x$, $t$ and $T$.
Suppose $r^2\ge T$, then we take $t=r^2$ in (2.20), we have
$$
\int_0^{r}sk_T(x,s)ds\le C_6\log (r+2)\tag 2.21
$$
for some constant $C_6$ independent of $x,t,T$.
Suppose $r^2\le T$, then by (2.19), we have
$$
\int_0^{r}sk_T(x,s)ds\le C_3\frac{r^2}{T+1}\le C_7\log (r+2)\tag2.22
$$
where $C_7$ is a  constant   independent of $x,t,T$.
(b) follows from (2.21) and (2.22).

To prove (a), if $2(1-\epsilon)/(1+\epsilon)<1$,
 the proof is similar to the proof of (b).
If $2(1-\epsilon)/(1+\epsilon)\ge1$,  the assumption in (a) implies
that $\int_0^\infty k(x,s)ds\le C_8$ for all $r$ and for all $x$.
By Theorem 1.3, for any $t$ we can solve the Poincar\'e-Lelong
equation $\i\ddbar u=\Ric(g(t))$ with $|\nabla u|(x,t)\le C_9$ for some
constant independent of $x$ and $t$. By Theorem 2.1 in \cite{N-S-T},
the result follows.
\enddemo

%\end

%Feb 28, 2002
\input amstex
\documentstyle{amsppt}
\magnification=1200
\hsize=13.8cm
\def \b {\beta}
\def\i{\sqrt{-1}}
\def\Ric{\text{Ric}}
\def\lf{\left}
\def\ri{\right}
\def\bbar{\bar \beta}
\def\a{\alpha}
\def\g{\gamma}
\def\p{\partial}
\def\delbar{\bar\delta}
\def\ddbar{\partial\bar\partial}

\def\C{\Bbb C}

\def\vp{\varphi}
\def\tD{\tilde \Delta}

\def\bb{{\bar\beta}}
\def\abb{{\alpha\bar\beta}}

\def\i{\sqrt {-1}}
\def\tD{\widetilde \Delta}
\def\tn{\widetilde \nabla}
\def\tg{\widetilde g}
\def \D {\Delta}

\def\aint{\frac{\ \ }{\ \ }{\hskip -0.4cm}\int}
\subheading{\S3 Asymptotic behavior of the volume element}

In \S2, we gave some estimates of the volume element $-F(x,t)$ under the K\"ahler-Ricci flow. In general, $-F(x,t)$ has no limit as $t\to\infty$ unless the original manifold is flat. In this section, we will use the Poincar\'e-Lelong equation and the results in \cite{N-S-T} to obtain information on asymptotic behavior of the rescaled volume element $-F(x,t)+F(x_0,t)$. Let us assume that $\lf(M^m,g_\abb(x)\ri)$ is a complete noncompact K\"ahler manifold with bounded nonnegative holomorphic bisectional curvature. As before, denote
$$
k(x,r)=\aint_{B(x,r)}\Cal R_0dV_0
$$
where $\Cal R_0$ is the scalar curvature of $g_\abb$. We also assume that
$$
k(x,r)\le k(r)\tag3.1
$$
for all $x\in M$, with $\int_0^\infty k(r)dr<\infty$. By Theorem 1.3,  (1.1) has a
long time solution $g_\abb(x,t)$ with $g_\abb(x,0)=g_\abb(x)$. On the other hand, by the result in \cite{N-S-T}, there is a unique function $u$ such that
$$
\i\ddbar u_0=\Ric(g(0))\tag3.2
$$
with $u_0(o)=0$ and $|u_0|=o(r)$. We have the following:
\proclaim{Theorem 3.1} Let $x_0\in M$ be a fixed point. For any $t_j\to\infty$, there is a subsequence, which is also denoted by $t_j$, such that
$$
\lim_{j\to\infty}\lf(F(x,t_j)-F(x_0,t_j)\ri)=u_0(x)-u_0(x_0)-v(x)
$$
where $u_0$ is the function in (3.2) and $v(x)$ is a pluriharmonic
function of at most linear growth (with respect to the initial metric). The convergence is uniform on compact sets. If in addition,
$\int_0^r sk(x,s)\, ds \le C(1+r)^{1-\epsilon}$ with $\epsilon>1/3$, then
$$
\lim_{t\to \infty}\lf(F(x,t)-F(x_0,t)\ri)=u_0(x)-u_0(x_0)
$$
and the convergence  is uniform on compact sets of $M$.
\endproclaim
\demo{Proof} Let $h(x,t)=\lf(u_0(x)-F(x,t)\ri)-\lf(u_0(x_0)-F(x_0,t)\ri)$. By Theorem 1.3, there exists a constant $C_1$ such that for all $(x,t)\in M\times[0,\infty)$
$$
\lf|\tn h(x,t)\ri|\le \lf|\nabla h(x,t)\ri|\le C_1,\tag3.3
$$
where $\tn h$ is the gradient with respect to the initial metric $g(0)$, and  we have used the fact that $g_\abb$ is nonincreasing. Since $h(x_0,t)=0$ for all $t$, it is easy to see that for any $t_j\to\infty$, there is a subsequence, which will be denoted by $t_j$ again, such that
$$
\lim_{j\to\infty}h(x,t_j)= v(x)
$$
for some Lipschitz continuous function $v(x)$ on $M$ with bounded gradient.
Since
$$
\tD h(x,t)=g^{\abb}(x,0)R_{\abb}(x,t)\tag3.4
$$
for all $x$, where $\tD$ is the Laplacian with respect to $g(0)$,
$$
0\le g^{\abb}(x,0)R_{\abb}(x,t)\le g^{\abb}(x,t)R_{\abb}(x,t)=\Cal R(x,t),
$$
Since by Corollary 2.4,  $\lim_{t\to\infty}\Cal R(x,t)=0$ uniformly on $M$, we   conclude that $v(x)$ is a
harmonic function of at most linear growth.
Notice that $h(x,t)$   is plurisubharmonic. Thus $v$ is also plurisubharmonic. Together with
the fact that it is also harmonic, $v$ must be pluriharmonic.

Suppose $\int_0^r sk(x,s)\, ds \le C(1+r)^{1-\epsilon}$ with $\epsilon>1/3$. Then by Theorem 2.3, we have
$$
\int_0^rsk_t(x,s)ds\le C_2 (1+r)^{\delta}\tag3.5
$$
for some constant $C_2>0$ independent of $x$ and $t$. Here $k_t(x,s)=\aint_{B_t(x,s)}\Cal RdV_t$ and $\delta=2(1-\epsilon)/(1+\epsilon)<1$. By Theorem 1.2
 in \cite{N-S-T} and the fact  that $h(x,t)=o\lf(r_t(x,x_0)\ri)$ for fixed $t$, we can conclude from (3.5) that
$$
h(x,t)\le C_3(1+r_t(x,x_0))^\delta\le C_3(1+r_0(x,x_0))^\delta
$$
for some constant independent of $t$. Hence the harmonic function $v(x)$ is of sublinear growth and must
 be constant by \cite{C-Y}. Since $v(x_0)=0$, $v$ must be identically zero.
\enddemo

In \cite{Sh2} and later in \cite{C-Z}, it was proved that if $M$ is a complete noncompact K\"ahler manifold with   {\it positive} and bounded holomorphic bisectional curvature such that the scalar curvature satisfies $\aint_{B(x,r)}\Cal R_0\le k(r)$ for all $x$ and $r$ with  with $k(r)\le C(1+r)^{-1-\epsilon}$, $\epsilon>1/2$, then the long time  solution of the K\"ahler-Ricci flow subconverges after rescaling in the following sense. Let $x_0$ be a fixed point in $M$ and let $v$ be a fixed   $(1,0)$ vector at $x_0$ with unit length with respect to the initial metric. Let $\widehat g_\abb(x,t)=g_\abb(x,t)/g_{v\bar v}(x_0,t)$. Then for any $t_j\to\infty$, we can find a subsequence, also denoted by $t_j$, such that $ \widehat g_\abb(x,t_j)$ converge uniformly on compact sets of $M$ to a flat K\"ahler metric. However, as pointed out in \cite{C-Z}, it is unclear whether the metric is complete. Using Theorem 3.1, we can get some preliminary estimates for the limiting metric.
\proclaim{Proposition 3.1} Let $(M^m,g_\abb)$ be a complete noncompact K\"ahler manifold with  positive and bounded holomorphic bisectional curvature such that the scalar curvature $\Cal R_0$ satisfies
$$
\aint_{B(x,r)}\Cal R_0dV_0\le k(r)
$$
for all $x$ and $r$, where $k(r)\le C(1+r)^{-1-\epsilon}$ with $\epsilon>1/2$. Let   $g_\abb(x,t)$ be the long time solution of (1.1) with $g_\abb(x,0)=g_\abb(x)$.
\roster
\item"{(a)}" The rescaled metrics
$$
\tg_\abb(x,t)=e^{-\frac{F(x_0,t)}m}g_\abb(x,t)
$$
subconverge to a flat K\"ahler metric $h_\abb$ on $M$. The convergence is uniform on compact sets, where $x_0$ is a fixed point and
$$
F(x,t)=\log \frac{\det(g_\abb(x,t))}{\det(g_\abb(x,0))}.
$$
\item"{(b)}" If, in addition,   $\epsilon=1$ and $\Cal R_0(x)\le Cr^{-2}_0(x)$, where $r_0(x)$ is the distance function from $x_0$ with respect to the initial metric, then
$$
\frac{\det(h_\abb(x,t))}{\det(g_\abb(x,0))}\ge C'r_0^a(x)-C''\tag3.6
$$
for some positive constants $a$, $C'$ and $C''$. In particular, the maximal eigenvalue $\lambda_{\text{\rm max}}(x)$ of $h_\abb(x)$ with respect to $g_\abb(x,0)$ satisfies
$$
\lambda_{\text{\rm max}}(x)\ge C'''r_0^{\frac am}(x)\tag3.7
$$
for some positive constant $C'''$, provided  $r_0(x)$ is large enough.
\endroster
\endproclaim
\demo{Proof} Part (a) follows from the results in
 \cite{Sh2, C-Z}.
Since
$$
\log\frac{\det(h_\abb(x,t))}{\det(g_\abb(x,0))}=\lim_{t\to\infty}\lf(F(x,t)-F(x_0,t)\ri),
$$
by Theorem 3.1 we have
$$
\log\frac{\det(h_\abb(x,t))}{\det(g_\abb(x,0))}=u(x)-u(x_0)\tag3.8
$$
where $u(x)$ is the solution for the Poincar\'e-Lelong equation obtained in \cite{N-S-T, Theorem 5.1}. Since $M$ is nonflat, by Remark 2.2 and Corollary 2.3,  we have
$$
\liminf_{r\to\infty}\frac{\int_0^rsk(x_0,s)ds}{\log r}>0.\tag3.9
$$
By \cite{N-S-T, Corollary 1.1}, (3.8) and (3.9), we  conclude that (3.6) is true.

(3.7) follows from (3.6) immediately.
\enddemo

%\end

\enddocument